\documentclass{amsart}
\theoremstyle{plain}
\newtheorem{Thm}{Theorem}
 
\newtheorem{Lem}[Thm]{Lemma}
\newtheorem{Prop}[Thm]{Proposition}

\numberwithin{equation}{section}

\newcommand{\CC}{{\mathbb{C}}}
\newcommand{\RR}{{\mathbb{R}}}
\newcommand{\QQ}{{\mathbb{Q}}}

\newcommand{\NN}{{\mathbb{N}}}
\begin{document}

\title[Reinhardt domains with a cusp at the origin]
{Reinhardt domains with a cusp at the origin}
\author{Oscar Lemmers}
\email[Oscar Lemmers]{lemmers@science.uva.nl}
\author{Jan Wiegerinck}
\email[Jan Wiegerinck]{janwieg@science.uva.nl}
\date{December 31, 2001}

\begin {abstract}
Let $\Omega$ be a bounded pseudoconvex Reinhardt domain in $\CC ^2$ with many 
strictly pseudoconvex points and logarithmic image $\omega$. It was known that
the maximal ideal in $H^{\infty}(\Omega)$ consisting of all functions 
vanishing at $(p_1,p_2) \in \Omega$ is generated by the coordinate functions 
$z_1 - p_1$, $z_2 - p_2$ (meaning that
one can solve the Gleason problem for $H^{\infty}(\Omega)$) if
$\omega$ is bounded. We show that one can solve Gleason's problem for 
$H^{\infty}(\Omega)$ as well if there are positive numbers $a$, $b$ and a
positive rational number $\frac{k}{l}$ such that $\Omega$ looks like
$\{(z_1,z_2) \in \CC^2 : a |z_2|^l \leq |z_1|^k \leq b |z_2|^l \}$ for small
$z$. 
\end{abstract}
\maketitle

\section{Introduction}

\noindent Let $\Omega$ be a bounded domain in $\CC^n$, let $p = (p_1,\ldots,
p_n)$ a point in $\Omega$. Recall the Gleason problem, cf. \cite{glea} : is
the maximal ideal in $A(\Omega)$ (or $H^{\infty}(\Omega)$) consisting of
functions vanishing at $p$ generated by the (translated) coordinate functions
$z_1 - p_1$, $\ldots$, $z_n - p_n$ ? We say that one can solve the Gleason
problem if this is indeed the case for every $p \in \Omega$. Gleason mentioned
the difficulty of solving this problem even for such a simple domain as the
unit ball $B(0,1)$ in $\CC^2$, $p=(0,0)$. This case was solved by Leibenzon
(\cite{khen}), who gave a solution tot the Gleason problem for every convex 
domain in $\CC^n$ with $C^2$ boundary.\\
Kerzman and Nagel (\cite{kena}) used 
sheaf-theoretic methods and  estimates on the solutions of $\overline
{\partial}$-problems to solve the Gleason problem for $A(\Omega)$, where 
$\Omega$ is a bounded strictly pseudoconvex domain in $\CC^2$ with 
$C^4$ 
boundary. Lieb (\cite{lieb}) independently solved the
Gleason problem for $A(\Omega)$ on bounded strictly pseudoconvex domains in 
$\CC^n$ with $C^5$ boundary;
\O vrelid improved this in \cite{oevr} to $C^{2}$ boundary. See also Henkin
(\cite{khen}) and Jakobczak (\cite{jako}).\\
In $\CC^2$ the Gleason problem was also solved for domains of finite 
type (\cite{foev}, \cite{noel}).
Backlund and 
F\"{a}llstr\"{o}m showed (\cite{bafa3}) that there exists an 
$H^{\infty}$-domain of holomorphy on which the Gleason problem is not solvable.
In \cite{bafa4} Backlund and F\"{a}llstr\"{o}m used ideas similar to those of 
Beatrous (\cite{beat}), to solve the Gleason problem for $A(\Omega)$ if 
$\Omega$ is a bounded pseudoconvex Reinhardt domain in $\CC^2$ with $C^2$
boundary that contains the origin. These ideas were expanded by the authors
(\cite{lewi}),
who solved the Gleason problem for both $A(\Omega)$ and $H^{\infty}(\Omega)$
if $\Omega$ is a bounded Reinhardt domain in $\CC^2$ with $C^2$ boundary. Thus 
the domain does not need to be pseudoconvex, and the condition that it contains
the origin could also be dropped. The condition of $C^2$ boundary could be 
weakened quite a lot, since it was only the behavior of the domain at 
the origin that was important. In this paper, we consider bounded pseudoconvex
Reinhardt domains $\Omega$ in $\CC^2$ that for small $z$ look like
\[\{(z_1,z_2) : a < \left|\frac{z_1^k}{z_2^l}\right|<b\}, \; k,l \in \NN^+, 
a,b \in \RR^+,\]
and are rounded of strictly pseudoconvexily. Thus, $\partial \Omega$ is 
non-smooth near the origin. We solve the Gleason problem for
$H^{\infty}(\Omega)$ in a way like \cite{lewi}. More detailed,
we divide the 
domain in two parts. On one part the problem is solved by splitting $f$ into 
functions for which an explicit solution is constructed. Adding these explicit 
solutions then gives a solution to the Gleason problem for $f$ on this part of
$\Omega$. On the other
part, the problem is solved using the $\overline{\partial}$-methods 
of \cite{lewi}. Then we patch the two local solutions together 
to a 
global solution, using a new $\overline{\partial}$-result.\\
We conclude by solving the Gleason problem for $H^{\infty}(\Omega)$ on the
Hartogs triangle and related domains.

\section{Definitions}

\noindent
We let $\CC^*$ stand for $\CC \setminus \{0\}$. 
Let \[L : (\CC^*)^2 \rightarrow \RR^2, \; L(z_1,z_2):=(\log |z_1|,\log|z_2|).\]
Throughout this paper $\Omega$ will be a bounded Reinhardt domain in 
$\CC^2$. 
We denote its logarithmic image $L(\Omega \cap (\CC^*)^2)$ by $\omega$. The 
boundary of $\Omega$
and $\omega$ will be denoted by $\partial \Omega$ and $\partial \omega$
respectively, while $S(\Omega)$ shall stand for the strictly
pseudoconvex boundary points of $\Omega$ that are $C^5$.\\
We denote the
derivative of a function $g$ with respect to the $j$'th coordinate with $D_jg$.
The interior and the closure of a set $V$ are denoted by $V^{\circ}$ and
$\overline{V}$ respectively. We denote the set \[\{(z_1,z_2) \in \CC^2 : 
g(z_1,z_2) = c\} \] by $[g(z_1,z_2)=c]$, and use a similar notation with e.g. 
$\leq$ instead of $=$.

\vskip5mm \noindent
{\bf Definition.} We say that $\Omega$ is an 
$A$-domain, if $\Omega$ is a bounded
pseudoconvex Reinhardt domain in $\CC^2$ such that 
\begin{itemize} 
\item There exist $a$, $b$, $\epsilon \in \RR^{+}$, $k$, $l \in \NN^{+}$, with
\[\Omega \cap 
B(0,\epsilon) = \{(z_1,z_2) \in B(0,\epsilon) : a < \left|\frac{z_1^{k}}
{z_2^l}\right| <b \}\]
\item The boundary points of $\Omega$ outside $\overline{B(0,\epsilon)}$ are 
all $C^5$ and strictly pseudoconvex.
\end{itemize}

\noindent{\bf Definition.} Let $U \subseteq \RR^n$ be an open set. For $0 <
\alpha < 1$ we define 
\[\Lambda_{\alpha}(U) = \{f \in C(U) : \sup_{x,x+h \in U} |f(x+h) - f(x)|
/|h|^{\alpha} + ||f||_{L^{\infty}(U)}  \] 
\[= ||f||_{\Lambda_{\alpha}}(U) < \infty\}.\]

\section{Solving a Cauchy-Riemann equation}

\noindent
The goal of this section is to prove the following theorem.

\begin{Thm}{\label{Thm:aholder}}
Let $\Omega$ be an $A$-domain. Suppose that $f$ is a 
$\overline{\partial}$-closed $(0,1)$-form with coefficients that are smooth and
bounded on $\Omega$, and that $\text {supp} f \cap \overline{B(0,
\epsilon)} = \emptyset$. Then there exists a $u \in \Lambda^{1/2}(\Omega)$ 
with $\overline{\partial}u = f$.
\end{Thm}

\noindent From this follows immediately that this $u$ is bounded on $\Omega$.
Note that under the assumptions of the theorem, the support of $f$ near the
boundary lies only near the strictly pseudoconvex points.
The setup of the
proof is very similar to that of the standard result on strictly
pseudoconvex domains with $C^5$ boundary. We will follow the book of Krantz 
(\cite{kran}), sections 5.2 and 9.1-9.3 (10.1-10.3 in the new edition). The
proof is subdivided in a series of lemmas. Proofs are given or indicated if
there is a difference with the standard situation, otherwise we 
refer to \cite{kran}. We do realize that the reader who is not that familiar 
with $\overline{\partial}$-problems will not be very happy about this decision.
In our opinion the alternative, copying over 25 pages word by word, would be
even worse.\\ 
Both in our case and the standard case, one has to
construct holomorphic support functions $\Phi(\cdot,P)$. Estimates on it
are derived by solving a $\overline{\partial}$-problem using the 
$L^2$-technique with weights of H\"{o}rmander (\cite{horm1}). In our case, we 
use that the 
$A$-domain $\Omega$ is contained in a slightly larger $A$-domain 
$\Omega_{1/n}$. The necessary estimate on a ball $B$ around the origin is
derived by a smart choice of the weight function $\phi$. The estimate on
$\Omega \setminus B$ is derived using that $\Omega \setminus B$  is compactly
contained in $\Omega_{1/n} \setminus B$. Compare this to the strictly
pseudoconvex case, where one uses that the domain is compactly contained in a 
strictly pseudoconvex domain that is strictly larger.

\vskip5mm \noindent We fix an $A$-domain $\Omega$.
Let $\epsilon$ be the smallest number such that $\partial \Omega \setminus 
\overline{B(0, \epsilon)}$ contains only 
strictly pseudoconvex points. We set $V:= \{w \in \partial \Omega : |w| > 
\epsilon\}$; then $V$ contains only strictly pseudoconvex points.
Let $\rho : \CC^2 \rightarrow \RR$ be a defining function for $\Omega$ that
is $C^5$ and strictly plurisubharmonic on a neighborhood of $V$.
The function $L : \CC^2 \times \CC^2 \setminus 
\overline{B(0,\epsilon)} \rightarrow \CC$ given by 
\[L_P(z)= L(z,P):= \rho(P) + \sum_{j=1}^{2} \frac{\partial \rho}{\partial z_j}
(P)(z_j - P_j) \]\[+ \frac{1}{2} \sum_{j,k=1}^{2} \frac{\partial ^2 \rho(P)}
{\partial z_j \partial z_k}(P)(z_j - P_j)(z_k - P_k)\]
is known as the 
Levi polynomial at $P$. It has the following properties :

\begin{enumerate}
\item For all $P \in \CC^2 \setminus \overline{B(0,\epsilon)}$, the function $z
\mapsto L(z,P)$ is holomorphic (it is even a polynomial).
\item For all $P \in V$, there is a neighborhood $U_P$ such that if $z \in
\overline{\Omega} \cap \{w \in U_P : L_P(w) = 0\}$ then $z=P$.
\end{enumerate}

\noindent
The goal is to construct for every $P \in V$ a holomorphic support function 
$\Phi(\cdot,P)$. This is a smooth function on $\Omega \times V$ that
is holomorphic in the first variable, such that $\Phi(z,P)=0 \Leftrightarrow
z=P$. Thus, this function should have the first property of the 
Levi polynomial at $P$. The difference is that one does not have to restrict 
in (2) to a small neighborhood of $P \in V$. The construction of these 
functions $\Phi(\cdot, P)$ will be done via some lemmas.

\vskip5mm \noindent
Choose $\gamma$, $\delta > 0$ such that 
\[\sum_{j,k=1}^{n} \frac{\partial ^2 \rho}{\partial z_j \partial 
\overline{z}_k}(P)v_j \overline{v}_k \geq \gamma |v|^2 \quad \forall 
P \in \{z \in \CC^n \setminus \overline{B(0,\epsilon)} : |\rho(z)| < \delta\}, 
v \in \CC^n.\]

\begin{Lem}{\label{Lem:lambda}}
There is a $\lambda > 0$ such that if $P \in V$ and $|z-P| < \lambda$,
then $2 \Re L_P(z) \leq \rho(z) - \frac{\gamma|z-P|^2}{2}$.
\end{Lem}

\noindent For every $n \in \NN$, we shall now define  $A$-domains 
 $\Omega_{1/n}$ that are close to $\Omega$. That is :
\[\Omega_{1/n} \cap B(0,\epsilon) = \{(z_1, z_2) \in B(0,\epsilon) : 
(1 - 1/n)a < \left|\frac {z_1^{k}}{z_2^{l}}\right| < (1 + 1/n)b\},\] and
$\Omega_{1/n}$ is rounded off strictly pseudoconvexily, having a $C^5$ 
defining 
function $\rho_{1/n}$ on a neighborhood $U$ of $\Omega_{1/n} \setminus 
B(0,\epsilon)$ such that
\begin{itemize}
\item $\Omega \subset \Omega_{1/n}$, $\partial \Omega \cap \partial 
\Omega_{1/n} = \{0\}$
\item $\Omega_{1/(n+1)} \subset \Omega_{1/n} \quad \forall n \in \NN$,
$\partial \Omega_{1/(n+1)} \cap \partial \Omega_{1/n} = \{0\}$
\item $\lim_{n \rightarrow \infty} ||\rho_{1/n} - \rho||_{C^5(U)} = 0.$
\end{itemize}
We also construct $A$-domains $\Omega_{-1/n}$ that are close to $\Omega$. That 
is :
\[\Omega_{-1/n} \cap B(0,\epsilon) = \{(z_1, z_2) \in B(0,\epsilon) : 
(1 + 1/n)a < \left|\frac {z_1^{k}}{z_2^l}\right| < (1 - 1/n)b\},\] and
$\Omega_{-1/n}$ is rounded off strictly pseudoconvexily, having a $C^5$ 
defining function $\rho_{-1/n}$ on a neighborhood $U$ of $\Omega \setminus 
B(0,\epsilon)$ such that
\begin{itemize}
\item $\Omega_{-1/n} \subset \Omega$, $\partial \Omega \cap \partial 
\Omega_{-1/n} = \{0\}$
\item $\Omega_{-1/n} \subset \Omega_{-1/(n+1)} \quad \forall n \in \NN$,
$\partial \Omega_{-1/(n+1)} \cap \partial \Omega_{-1/n} = \{0\}.$
\item $\lim_{n \rightarrow \infty} ||\rho_{-1/n} - \rho||_{C^5(U)} = 0.$
\end{itemize}
This is possible, cf. the setup in \cite{lewi} : we only need to
consider convex domains in $\RR^2$ instead of pseudoconvex Reinhardt domains
in $\CC^2$.\\
We choose $n \in \NN$ such that $||\rho_{1/n} - \rho||_{C^5(U)} \leq 
\frac{\gamma \lambda^2}{20}$ (where $\lambda$ is the constant of lemma 
\ref{Lem:lambda}). We may assume that $||\rho_{1/n} - \rho||_{C^5(U)} < 
\lambda < \delta <1$.

\begin{Lem}
If $P \in V$, $z \in \Omega_{1/n}$, $\lambda/3 \leq |z-P| \leq 2\lambda/3$,
then $\Re L_P(z) < 0$.
\end{Lem}

\noindent
Let $\eta : \RR \rightarrow [0,1]$ be a $C^{\infty}$ function that satisfies
$\eta(x)=1$ for $x \leq \lambda/3$, $\eta(x)=0$ for $x \geq 2\lambda/3$.

\begin{Lem}
Let $P \in V$. The $(0,1)$-form 
\[f_P(z) = \left\{ \begin{array}{ll} -\overline{\partial}_z(\eta(|z-P|)) \cdot
\log L_P(z) & \mbox{\!\!,\quad if $|z-P| < \lambda, z \in \Omega_{1/n}$}\\
0 & \mbox{\!\!,\quad if $|z-P| \geq \lambda, z \in \Omega_{1/n}$}\\
\end{array}\right.\]
is well defined (if we take the principal branch for the logarithm) and has
$C^{\infty}$ coefficients for $z \in \Omega_{1/n}$. Furthermore, 
$\overline{\partial}_z f_P(z) = 0$ on $\Omega_{1/n}$.
\end{Lem}

\begin{Lem}{\label{Lem:dbaradom}}
Let $f$ be a $\overline
{\partial}$-closed $(0,1)$-form on $\Omega_{1/n}$ with $C^1$ coefficients
that are bounded.
Suppose that $\overline{B(0, \epsilon)} \cap \text{supp} f = 
\emptyset$. Then there exist a $C_{\epsilon}$ (that does not 
depend on $f$) and a function $u$ with $\overline{\partial}u = f$
such that \[||u||_{L^{\infty}(\Omega_{1/2n})} \leq C_{\epsilon}||f||_
{L^{\infty}(\Omega_{1/n})}.\]
\end{Lem}

\begin{proof}
If $f$ is identically zero, we are done. So assume that 
$||f||_{L^{\infty}(\Omega_{1/n})} > 0$. 
We choose a weight $\phi$ that blows up near the boundary of $\Omega_{1/n}$. 
Then
we add several times $\log |z|$ such that $e^{- \phi(z)}$ will behave like
$|z|^{-k}$ (this $k$ will be chosen later). We let $u$ be the solution of the
$\overline{\partial}$-equation on $\Omega_{1/n}$ for the weight $\phi$, as
constructed by H\"{o}rmander (\cite{horm1}). Then 
\[\int_{\Omega_{1/n}} \frac{|u(z)|^2 e^{- \phi(z)}}{(1 + |z|^2)^2} d\lambda
\leq \int_{\Omega_{1/n}} |f(z)|^2 e^{- \phi(z)} d\lambda < \infty.\]
The first inequality is the estimate of H\"{o}rmander, the second one holds
because $f$ has bounded coefficients. We start by showing that the 
assumption 
that there is a sequence $\{z_n\}^{\infty}_{n=1}$ in $\Omega_{1/(2n)}$ that 
converges to
$0$ such that $|u(z_n)| \geq ||f||_{L^{\infty}(\Omega_{1/n})}$ leads to a
contradiction. This yields an estimate for $||f||_{\infty}$ near the origin.\\
There are constants $R$, $\beta > 0$ such that 
\[z \in \Omega_{1/(2n)} \cap B(0,\epsilon) \Rightarrow B(z,R(|z|^{\beta}) 
\subset \Omega_{1/n}.\]
Thus for large $n$, one has that $B(z_n, R|z_n|^{\beta})$ is contained
completely in $B(0,\epsilon) \cap \Omega_{1/n}$.
We choose $k > 4 \beta$.
We assumed that $f$ has no support on $\overline{B(0,\epsilon)} \cap 
\Omega_{1/n}$, thus $u$ is holomorphic there. We now
apply the mean value inequality on $u$.
\[\infty > \int_{\Omega_{1/n}} |f(z)|^2 e^{- \phi(z)} d\lambda \geq 
\int_{\Omega_{1/n}} 
\frac{|u(z)|^2 e^{- \phi(z)}}{(1 + |z|^2)^2} d\lambda \]\[
> \int_{B(z_n, R|z_n|^
{\beta})} \frac{|u(z)|^2 e^{- \phi(z)}}{(1 + |z|^2)^2} d\lambda > C
\frac{|u(z_n)|^2 R^4 |z_n|^{4 \beta}}{|z_n|^k} > C' |z_n|^{4 \beta - k} 
\rightarrow \infty\]
if $n \rightarrow \infty$. Thus there is a $\delta$ with $0 < \delta < 
\epsilon$ such that $|u(z)| \leq 
 ||f||_{L^{\infty}(\Omega_{1/n})}$ for $z \in \Omega_{1/(2n)} \cap 
B(0,\delta)$.\\
Now we shall make the appropriate estimate on $\Omega_{1/(2n)} \setminus
B(0,\delta)$.
Remember the H\"{o}rmander construction (\cite{kran}, chapter 4), with
$\phi$, $\phi_1$, $\phi_2$ and
\[T = \overline{\partial}_{0,0} : L^2_{(0,0)}(\Omega_{1/n}, \phi_1) \rightarrow
L^2_{(0,1)}(\Omega_{1/n}, \phi_2).\]
Then
\[\sup _{\Omega_{1/(2n)} \setminus B(0, \delta)} |u| \leq C(||u||_{L^2(
\Omega_{3/(4n)} 
\setminus B(0, \delta))} + ||\overline{\partial} u||_{L^{\infty}_{(0,1)}
(\Omega_{3/(4n)} \setminus B(0, \delta))})\]
\[\leq C' (||u||_{L^2(\Omega_{1/n} \setminus B(0, \delta), \phi_1)} + 
||f||_{L^{\infty}_
{(0,1)} (\Omega_{3/(4n)} \setminus B(0, \delta))})\]
\[\leq C'' (||f||_{L^2_{(0,1)}(\Omega_{1/n} \setminus B(0, \delta), \phi_2)} + 
||f||_{L^{\infty}_{(0,1)} (\Omega_{3/(4n)} \setminus B(0, \delta))})\]
\[\leq C''' ||f||_{L^{\infty}_{(0,1)}(\Omega_{1/n} \setminus B(0, \delta))}\]
since $e^{- \phi_2(z)}$ tends to zero as $z$ tends to a boundary point of
$\Omega_{1/n}$ that is non-zero.
\end{proof}

\noindent
For every $P \in V$, we let $u_P$ be a solution of $\overline{\partial}u_P = 
f_P$ that satisfies the estimate above. We now define
\[\Phi(z,P) = \left\{ \begin{array}{ll} L_P(z) \cdot \exp(u_P(z)) & \mbox{\!\!,
\quad if $|z-P| < \lambda/3$} \\
\exp(u_P(z) + \eta(|z-P|)\log L_P(z)) & \mbox{\!\!, \quad if $\lambda/3 \leq
|z-P| < \lambda$} \\
\exp(u_P(z)) & \mbox{\!\!, \quad if $\lambda \geq |z-P|$}\\ \end{array} 
\right.\]
We proceed to show that these functions $\Phi(\cdot,P)$ are holomorphic support
functions.

\begin{Lem}{\label{Lem:holsup}}
For every $P \in V$,
the function $\Phi(\cdot,P)$ is holomorphic on $\Omega_{1/n}$. For fixed $z \in
\Omega_{1/(2n)}$, $\Phi(z,\cdot)$ is continuous in $P$. There is a
$C > 0$, independent of $P$, such that for all $z \in \Omega_{1/(2n)}$ we have
\[\text{if} \quad |z-P| < \lambda/3, \quad \text{then} \quad |\Phi(z,P)| 
\geq C|L_P(z)|,\]
\[\text{if} \quad |z-P| \geq \lambda/3, \quad \text{then} \quad |\Phi(z,P)| 
\geq C.\]
\end{Lem}

\begin{proof}
The function $f_P$ is bounded on $\Omega_{1/(2n)}$ uniformly in $P$, hence 
$u_P$ is bounded on $\Omega_{1/(2n)}$
uniformly in $P$. Thus there is a $C > 0$ such that $|\exp
u_P(z)| \geq C$. Working this out yields the appropriate estimates.
\end{proof}

\begin{Lem}{\label{Lem:holsupbdd}}
For every $P \in V$ there exist functions $\Phi_1(z,P)$, $\Phi_2(z,P)$ that
are holomorphic in $z \in \Omega_{1/n}$ and a constant $C$ that does not
depend on $P$, such that
\[\Phi(z,P) = \Phi_1(z,P) (z_1 - P_1) + \Phi_2(z,P) (z_2 - P_2) \quad \forall
z \in \Omega_{1/n},\]
\[|\Phi_j(z,P)| \leq C \quad \forall z \in \Omega_{1/(2n)}, P \in V, j = 1,2.\]
\end{Lem}

\begin{proof}
We will follow the approach of Backlund and F\"{a}llstr\"{o}m in \cite{bafa4}.
A line with positive rational slope $\frac{k}{l}$ in $\RR^2$ passing through
$L(p)$ can be seen as the logarithmic image of the zero set of the polynomial 
$z_1^k p_2^l - z_2^l p_1^k$, while a line with negative rational slope 
$\frac{-k}{l}$ in $\RR^2$ passing through $L(p)$ can be seen as the logarithmic
image of the zero set of the polynomial $z_1^k z_2^l - p_1^k p_2^l$.\\
Fix $P \in V$.
We choose polynomials $g$ and $h$ such that $L(Z_g)$ and $L(Z_h)$ are lines
in $\RR^2$ that intersect
the boundary of $\Omega$ only in $V$, and $[g=0] \cap [h=0] \cap 
\Omega_{1/n} = \{P\}$.
Now choose a ball $U_0$ around $P$ that lies compactly in $\Omega_{1/n}$, and
choose open sets $U_1$, $U_2$ such that
\begin{itemize}
\item $\overline{\Omega_{1/n}} \subset \cup_i U_i$
\item For a certain positive number $\mu$ one has that $|g| > \mu$ on $U_1$,
$|h| > \mu$ on $U_2$.
\item $\overline{U_1} \cap \overline{U_2} \cap \overline{B(0,\epsilon)} =
\emptyset$.
\end{itemize}

\noindent Now choose functions $\phi_k \in C^{\infty}_0(U_k) \quad (k=0,1,2)$ 
such that $0 \leq \phi_k \leq 1$ and $\sum_{k=0}^2 \phi_k = 1$ on $\overline
{\Omega_{1/n}}$.
Recall that $\Phi(\cdot,P)$ vanishes at $z=P$. Because $\Phi(\cdot,P)$ is
holomorphic on $\Omega_{1/n}$, and $U_0 \subset \subset \Omega_{1/n}$, the
lemma of Oka-Hefer implies that there exist functions $\Phi_1^0(\cdot,P)$, 
$\Phi_2^0(\cdot,P) \in H^{\infty}(U_0)$ such that 
\[\Phi(z,P) = \Phi_1^0(z,P)(z_1 - P_1) + \Phi_2^0(z,P)(z_2 - P_2) \quad \quad
\forall z \in U_0.\]
We set \[\tilde{\Phi}_1^1(z,P):= \frac{\Phi(z,P)}{g(z)}, \tilde{\Phi}_2^1(z,P)
:=0,\] \[\tilde{\Phi}_1^2(z,P):=0, \tilde{\Phi}_2^2(z,P):=\frac{\Phi(z,P)}
{h(z)}.\]
Then $\tilde{\Phi}_j^i \in H^{\infty}(U_i \cap \Omega_{1/n})$ and
\[\Phi(z,P) = \tilde{\Phi}_1^i(z,P)g(z) + \tilde{\Phi}_2^i(z,P)h(z) \quad 
\forall z \in U_i, i \in \{1,2\} \quad \quad \quad(*).\]
Since $g$ is an analytic polynomial, vanishing at $P$, there are polynomials
$g_{1}$, $g_{2} \in H(\CC^{2})$ such that
$g = g_{1} (z_1 - P_1) + g_{2} (z_2 - P_2)$ on $\CC^{2}$. A similar
formula holds for $h$. Substituting this in $(*)$, we obtain the existence of 
functions $\Phi^{i}_{j} \in H^{\infty}(U_i \cap \Omega_{1/n})$, $i = 1$, $2$, 
such that
\[ \Phi(z,P) = \Phi_{1}^{i}(z,P) (z_1 - p_1) + \Phi_{2}^{i}(z,P) (z_2 - p_2) 
\; \text { on } \; \overline{U_i} \cap \Omega_{1/n}, \; \; \; i = 1, 2.\]
Therefore 
\[j_1 := \sum_{k=0}^{2} \phi_k \Phi_1 ^k \; \text { and } \; 
j_2 := \sum_{k=0}^{2} \phi_k \Phi_2 ^k \]
\noindent give a smooth solution of our problem. We want to find $u$ such that
\[\Phi_1 = j_1 + u(z_2 - P_2) \; \text { and } \; \Phi_2 = j_2 - u(z_1 - P_1) 
\quad \quad \quad (**)\] are in $H(\Omega_{1/n}) \cap L^{\infty}
(\Omega_{1/(2n)})$.
Define a form $\lambda$ as follows :
\[ \lambda := \frac {- \overline{\partial}j_1}{z_2 - P_2} = 
\frac {\overline{\partial}j_2}{z_1 - P_1} \]
\noindent This form $\lambda$ is a bounded $\bar{\partial}$-closed 
$(0,1)$-form on $\Omega_{1/n}$.
The support of $\lambda$ is contained in $\overline{U_i} \cap \overline{U_j}$, 
$i \not = j$. These sets all lie outside $\overline{B(0,\epsilon)}$.
Lemma \ref{Lem:dbaradom} gives the existence of a 
function $u \in L^{\infty}(\Omega_{1/(2n)})$ 
such that $\bar{\partial}u = \lambda$. 
With this $u$, $\Phi_1$, $\Phi_2$ as defined at $(**)$, 
 \[\Phi = \Phi_1 (z_1 - P_1) + \Phi_2 (z_2 - P_2)\] on $\Omega_{1/n}$
, and $\Phi_1(\cdot,P)$, $\Phi_2(\cdot,P)$ both belong to 
$H(\Omega_{1/n}) \cap L^{\infty}(\Omega_{1/(2n)})$.\\
For fixed $z \in \Omega_{1/(2n)}$, the function $\Phi(z,\cdot)$ depends 
continuously on $P$. Studying the 
construction above carefully, we see that we can choose $\Phi_1(z,\cdot)$ and 
$\Phi_2(z,\cdot)$ continuously in $P$ as well. Thus, because supp $\lambda \cap
\partial \Omega$ is compact,
there exists a uniform bound on $||\Phi_i(z,P)||_{\Omega_{1/(2n)}}$.
\end{proof}

\begin{Thm}
Let $\Omega$ be an $A$-domain. Let $f$ be a $\overline{\partial}$-closed 
$(0,1)$-form on an $A$-domain that contains $\overline{\Omega} \setminus \{0\}$
with $C^1$ coefficients. 
Suppose that $\text{supp} f \cap \overline{B(0,\epsilon)} = \emptyset$.
Then there is a function $u$ such that $\overline{\partial} u =f$, and 
\[||u||_{L^{\infty}(\Omega)} \leq C||f||_{L^{\infty}(\Omega)}.\]
\end{Thm}

\begin{proof}
Let 
$H_{\Omega}(f)(z)$ be the Khenkin solution to the $\overline{\partial}$ 
equation; then $\overline{\partial} H_{\Omega}(f) = f$.
To prove the necessary estimates, we start by
writing $f=f_1 d \overline{\zeta}_1 + f_2 d \overline{\zeta}_2$.
Then the Khenkin solution can be rewritten to 
\[H_{\Omega}(f)(z)= \frac{1}{4 \pi^2} \{\int_{\Omega} \frac{f_1(\zeta)\cdot
(\overline{\zeta}_1 - \overline{z}_1) + f_2(\zeta) \cdot (\overline{\zeta}_2 - 
\overline{z}_2)}{|\zeta - z|^4} \times d\overline{\zeta}_1 \wedge d\overline
{\zeta}_2 \wedge d \zeta_1 \wedge d \zeta_2 \]
\[- \int_{\partial \Omega} \frac{\Phi_1(z,\zeta)(\overline{\zeta}_2 - 
\overline{z}_2) - \Phi_2(z,\zeta)(\overline{\zeta}_1 - \overline{z}_1)}
{\Phi(z, \zeta)|\zeta - z|^2} \times (f_1(\zeta)d\overline{\zeta}_1 + f_2(
\zeta)d \overline{\zeta}_2) \wedge d \zeta_1 \wedge \zeta_2 \}\]
\[= \int_{\Omega} f_1(\zeta)K_1(z, \zeta) dV(\zeta) + 
\int_{\Omega} f_2(\zeta)K_2(z, \zeta) dV(\zeta) \]
\[+ \int_{\partial \Omega} f_1(\zeta)L_1(z, \zeta) dV(\zeta) + 
\int_{\partial \Omega} f_2(\zeta)L_2(z, \zeta) dV(\zeta)\]
where the identity defines the kernels. Now let $T$ be so large that $\Omega
\subseteq B(z,T)$ for every $z \in \Omega$. Then
\[\int_{\Omega} |K_i(z, \zeta|dV(\zeta) \leq \int_{B(z,T)} |z - \zeta|^{-3}
dV(\zeta) = C \int_{0}^{T} r^{-3}r^3 dr = C' \quad j=1,2.\]
Because $f$ has no support on $\overline{B(0,\epsilon)}$, one has that
\[\int_{\partial \Omega} f_j(\zeta)L_j(z, \zeta) dV(\zeta) = 
\int_{V} f_j(\zeta)L_j(z, \zeta) dV(\zeta) \quad j=1,2 .\]
Using lemmas \ref{Lem:holsup} and \ref{Lem:holsupbdd}, one can prove that 
\[\int_{V} |L_j(z, \zeta)| dV(\zeta) \leq D_j \quad j=1,2,\]
where the bounds are independent of $z \in \Omega$. 
This implies that 
\[||H_{\Omega}(f)||_{L^{\infty}(\Omega)} \leq (2 C' + D_1 + D_2) 
||f||_{L^{\infty}(\Omega)}.\]

\noindent Keeping in mind that $|\Phi_1(z,\zeta)|$ and $|\Phi_2(z,\zeta)|$
are bounded on $\Omega_{1/(2n)}$ uniformly in $\zeta$, one can simply follow 
\cite{kran}. 
\end{proof}

\noindent Repeating all the arguments used over there exactly, yields :

\begin{Thm}{\label{Thm:holder}}
Let $\Omega$ be an $A$-domain. Let $f$ be a $\overline{\partial}$-closed 
$(0,1)$-form on an $A$-domain that contains $\overline{\Omega} \setminus \{0\}$
with $C^1$ coefficients. 
Suppose that $\text {supp} f \cap \overline{B(0,\epsilon)} = 
\emptyset$. Then $H_{\Omega}(f)$ is well defined, continuous on $\overline
{\Omega}$ and
\[||H_{\Omega}(f)||_{\Lambda_{1/2}(\Omega)} \leq C||f||_{L^{\infty}(\Omega)}.\]
\end{Thm}

\begin{Thm}
Let $\Omega$ be an $A$-domain. Then there is an $N \in 
\NN$ such that if $n \geq N$, then
theorem \ref{Thm:holder} holds on $\Omega_{-1/n}$ with $C_{\Omega_{-1/n}} 
\leq 2 C_{\Omega}$.
\end{Thm}

\noindent
Now we give the proof of theorem \ref{Thm:aholder}. \begin{proof}
Let $\Omega$ be an $A$-domain. For $n \in \NN$ large, 
the stability result will apply on $\Omega_{-1/n}$. Now let $f$ be a $\overline
{\partial}$-closed $(0,1)$-form defined on $\Omega$ (not necessarily on a
neighborhood of $\overline{\Omega}$) with bounded $C^1$ coefficients. For
each sufficiently small $-1/n < 0$, the form $f$ satisfies the hypotheses
of theorem \ref{Thm:holder} on $\Omega_{-1/n}$. Therefore $H_{\Omega_
{-1/n}}(f)$ is well defined and satisfies $\overline{\partial}H_{\Omega_
{-1/n}}(f) = f$ on $\Omega_{-1/n}$. Moreover,
\[||H_{\Omega_{-1/n}}(f)||_{\Lambda_{1/2}(\Omega)} \leq C_{\Omega_
{-1/n}}||f||_{L^{\infty}(\Omega)} \leq 2 C_{\Omega}||f||_{L^{\infty}(\Omega)}
.\]
Thus, given a compact subset $K$ of $\Omega$, the functions 
$\{H_{\Omega_{-1/n}}(f)\}$ form an equicontinuous family on $K$ if $n$ is 
large. Of course, this family is also equi-bounded. By the Arzel\`{a}-Ascoli
theorem and diagonalization, we see that there is a
subsequence $H_{\Omega_{-1/j}}(f)$, $j=1$, $2$, $\ldots$, such that $H_{\Omega_
{-1/j}}(f)$ converges uniformly on compacta to a $u \in \Lambda_{1/2}
(\Omega)$ with $\overline{\partial}u = f$ on $\Omega$.
\end{proof}

\noindent {\bf Remark.} Note that theorem \ref{Thm:aholder} also holds for e.g.
a Reinhardt domain $\Omega$ that for small $z$ looks like 
\[\{(z_1,z_2) : 0 < |z_1^k| < |z_2^l|\},\]
and is rounded off strictly pseudoconvexily.

\section{Auxiliary results}

\begin{Lem}{\label{Lem:solfrac}}
Let $\Omega$ be a domain in $\CC^2$, let $(p_1,p_2) \in \Omega$, let $k$, $l 
\in \NN^*$.
Suppose that $\frac{z_1^k}{z_2^l} \in H^{\infty}(\Omega)$.
Let \[R_1(z_1,z_2) := \frac{1}{p_2^l} \frac{z_1^k - p_1^k}{z_1 - p_1},\] 
\[R_2(z_1,z_2) :=  \frac{1}{p_2^l} \frac{z_1^k}{z_2^l} \frac{p_2^l - z_2^l}
{z_2-p_2} .\]
Then \[\frac{z_1^k}{z_2^l} - \frac{p_1^k}{p_2^l} = R_1(z_1,z_2) (z_1-p_1) + 
R_2(z_1,z_2) (z_2-p_2),\] and $R_1$, $R_2 \in H^{\infty}(\Omega)$.
\end{Lem}

\begin{proof}
This can be checked by hand. 
\end{proof}

\begin{Lem}{\label{Lem:solpol}}
Let $P$ be a polynomial in $z_1$ and $z_2$ that vanishes at $(p_1,p_2) \in 
\CC^2$. There exist polynomials $P_1$, $P_2$ such that \[P(z_1,z_2) = P_1(z_1,
z_2)(z_1-p_1) + P_2(z_1,z_2)(z_2-p_2).\] 
\end{Lem}

\begin{proof}
For $(p_1,p_2)=(0,0)$, this follows immediately. For other points apply the 
appropriate coordinate transform.
\end{proof}

\begin{Lem}{\label{Lem:driehoek}}
Suppose there are points $t$, $u$, $v \in \partial \omega$ having 
neighborhoods $T$, $U$, $V \subset \partial \omega$
consisting only of strictly convex points of $\partial \omega$ respectively,
such that $L(p) \in Co(tuv)$. Then one can solve the Gleason problem for 
$H^{\infty}(\Omega)$ at $p$.
\end{Lem}

\begin{proof}
We choose, just as in lemma \ref{Lem:holsupbdd}, analytic polynomials $g$, 
$h$, open 
sets $U_0$, $U_1$, $U_2$ and a constant $\mu > 0$ such that:
\begin{itemize}
\item $[g(z) = 0] \cap [h(z) = 0] \cap \overline{\Omega} = \{p\}$ 
\item $U_0$ is strictly pseudoconvex, and $p \in U_0 \subset \subset 
\Omega$
\item $|g| > \mu$ on $U_1$, $|h| > \mu$ on $U_2$
\item $\overline{\Omega} \subset \cup _i U_i$
\item $U_i \cap U_j \cap B(0,\epsilon) = \emptyset, j=0,1,2.$ 
\end{itemize}
Now formulate the corresponding $\overline{\partial}$-problem, again as in 
lemma \ref{Lem:holsupbdd}. This yields a
bounded $(0,1)$-form that has only support outside $B(0, \epsilon)$. 
Applying theorem \ref{Thm:aholder} yields a bounded solution to the
$\overline{\partial}$-problem, and this can be used to solve the Gleason
problem in the standard way.
\end{proof}

\section{Dividing $\Omega$ in two pieces}{\label{section:pieces}}

\noindent Suppose that $\Omega$ is an $A$-domain, and that $p \in \Omega$.
Then the line with slope $\frac{k}{l}$ through $L(p)$ intersects $\partial
\omega$ in only one point $A$. This point is strictly convex. Thus there is a
line $N$ in $\RR^2$ with rational slope $\not=\frac{k}{l}$ that 
intersects $\partial \omega$ only at strictly convex points such that $A$ and
the part of $\omega$ in the third quadrant lie on different sides of $N$.
Say $N$ is given by the equation 
$y=\frac{-m}{n} x + r$, where $m$, $n \in \NN$. Then $N$ is the logarithmic
image of $[z_1^m z_2^n = e^{rn}]$. There is a $\delta > 0$ such that 
\[\{(x,y) \in \partial \omega, \tilde{r} \in [r - \delta, r] : y = 
\frac{-m}{n}x + \tilde{r} \} \subset S(\omega).\]
Let
\[\omega_1 := \{(x,y) \in \omega : y \geq \frac{-m}{n}x + r - \delta\},\] 
\[\omega_2 := \{(x,y) \in \omega : y \leq \frac{-m}{n}x + r\},\] 
and $\Omega_1$, $\Omega_2$ be $(\overline{L^{-1}(\omega_1)})^{\circ}$, 
$(\overline{L^{-1}(\omega_2)})^{\circ}$ respectively.
If $p$ lies in $\Omega_1$, everything is easy : apply lemma 
\ref{Lem:driehoek} to solve the Gleason problem for $H^{\infty}(\Omega)$ at 
$p$. 
\\In the rest of the article we shall assume that $p$ does not lie in 
$\Omega_1$. We will use that there is an
$\nu > 0$ such that $|z_1^m z_2^n - p_1^m p_2^n| > \nu$ for $(z_1,z_2) \in 
\Omega_1$ to obtain a local solution on $\Omega_1$. The next section
consists of the construction of a local solution on $\Omega_2$. Afterwards, 
the two local solutions will be patched together using the standard arguments.

\section{Constructing a local solution}{\label{section:constrloc}}

\noindent We fix $p= (p_1,p_2) \in \Omega_2$ and $f \in H^{\infty}(\Omega)$ 
that vanishes at $p$.
The main idea of the following construction is to project 
$(z_1, z_2)$ on the zero set of $\frac{z_1^k}{z_2^l} - \frac{p_1^k}{p_2^l}$, 
because
\[\frac{f(z_1, z_2) - f((\frac{p_1^k z_2^l}{p_2^l})^{1/k}, z_2)}{\frac{z_1^k}
{z_2^l} - \frac{p_1^k}{p_2^l}} (\frac{z_1^k}{z_2^l} - \frac{p_1^k}{p_2^l}) +
\frac{f((\frac{p_1^k z_2^l}{p_2^l})^{1/k}, z_2)}{z_2-p_2} (z_2-p_2)\]
comes close to being a solution for the Gleason problem. However, as there 
appear roots in the argument of the function, we lose in general the
holomorphy. We decompose $f$ in functions where one can take the
appropriate root. Then we solve the Gleason problem for those functions, add 
all these solutions and end up with a solution of the Gleason problem for $f$. 

\vskip5mm \noindent By 
$\zeta$ we denote the $(kn +lm)$'th elementary root of unity.

\begin{Lem}{\label{Lem:symmetry}}
Suppose $f$ is a bounded holomorphic function on $\Omega_2$. Then for every
$0 \leq i$, $j \leq kn + lm - 1$ there exist functions $f_{i,j} \in 
H(\Omega_2)$ such that :
\begin{itemize}
\item $z_1^i z_2^j f_{i,j}$ is bounded for $0 \leq i,j \leq kn +lm -1$
\item $f_{i,j}(z_1,z_2) = f_{i,j}(\zeta z_1,z_2) = f_{i,j}(z_1, \zeta z_2)$ 
for all $(z_1,z_2) \in \Omega_2$, $0 \leq i,j \leq kn + lm -1$
\item $f(z_1,z_2) = \sum_{i,j=0}^{kn + lm -1} z_1^i z_2^j f_{i,j}(z_1,z_2)$ 
for all $(z_1,z_2) \in \Omega_2$.
\end{itemize}
\end{Lem}

\begin{proof}
Let 
\[f_{i,j}(z_1,z_2) := \frac{1}{(kn + lm)^2 z_1^i z_2^j} \sum_{s,t=1}^
{kn + lm} \zeta^{-is -jt} f(\zeta^{s}z_1,\zeta^{t}z_2).\] 
The domain $\Omega$ does not contain points with a zero coordinate, hence
$f_{i,j}$ is well defined. Since $f$ is bounded, we see immediately that 
$z_1^i z_2^j f_{i,j}(z_1,z_2)$ is bounded as well. 
\[ (kn + lm)^2 f_{i,j}(\zeta z_1, z_2) = \frac{1}{(\zeta z_1)^i z_2^j} 
\sum_{s,t=1}^{kn +lm} \zeta^{-is-jt} f(\zeta^{s+1}z_1, \zeta^{t}z_2) = \]
\[\frac{1}{(\zeta z_1)^i z_2^j} \sum_{t=1}^{kn + lm} \zeta^{-jt} 
\sum_{s=2}^{kn +lm +1} \zeta^{-i(s-1)} f(\zeta^{s}z_1, \zeta^{t}z_2) = \] \[
\frac{\zeta^{i}}{(\zeta z_1)^i z_2^j} \sum_{t=1}^{kn + lm} 
\zeta^{-jt}\left( \zeta^{-i(kn + lm +1)} f(\zeta^{kn + lm +1}z_1,\zeta^{t} z_2)
+ \sum_{s=2}^{kn + lm} \zeta^{-is} f(\zeta^{s}z_1, \zeta^{t}z_2) \right) =\] 
\[=\frac{1}{ z_1^i z_2^j} \sum_{s,t=1}^{kn + lm} \zeta^{-is-jt} 
f(\zeta^{s}z_1, \zeta^{t}z_2) = (kn +lm)^2 f_{i,j}(z_1,z_2).\]
The equality $f_{i,j}(z_1, \zeta z_2) = f_{i,j}(z_1,z_2)$ can be proven 
similarly. 
Since \[\sum_{i,j=0}^{kn + lm -1} \zeta^{-is -jt} = \sum_{i=0}^{kn + lm -1} 
\zeta^{-is} \sum_{j=0}^{kn + lm -1} \zeta^{-jt} = \left\{ \begin{array}{ll} 
0 & \mbox {\quad $s,t \neq kn + lm$} \\ (kn +lm)^2 & \mbox{\quad $s,t=kl +mn$} 
\\ \end{array} \right. \]
we have that \[\sum_{i,j=0}^{kn +lm -1} z_1^i z_2^j f_{i,j}(z_1,z_2) = \frac{1}
{(kn +lm)^2} \sum_{i,j=0}^{kn +lm -1}
\sum_{s,t=1}^{kn + lm} \zeta^{-is -it} f(\zeta^{s}z_1, \zeta^{t}z_2) = \] 
\[=\frac{1}{(kn +lm)^2} 
\sum_{s,t=1}^{kn +lm} f(\zeta^{s}z_1, \zeta^{t}z_2) \sum_{i,j=0}^{kn +lm -1} 
\zeta^{-is -it} = f(z_1,z_2).\]
\end{proof}

\noindent
{\bf Remark.} There is a polynomial $P$ such that \[P(\zeta^s p_1, \zeta^t p_2)
= f(\zeta^s p_1, \zeta^t p_2) \quad \forall 1 \leq s,t \leq kn + lm.\] 
From lemma \ref{Lem:solpol} it follows that one can solve the Gleason problem 
for the function 
$f$ if and only if one can solve the Gleason problem for $f - P$.
The corresponding functions $(f - P)_{i,j}$ all vanish
at $p$. Hence we may assume from now on that $f_{i,j}$ vanishes at 
$p$. 

\begin{Lem}
The multi valued map $\pi$ given below, maps a point $(z_1,z_2) \in \Omega_2$
to the set $[\frac{z_1^k}{z_2^l}=\frac{p_1^k}{p_2^l}] \cap \Omega$. 
The function $f_{i,j} \circ \pi$ is a holomorphic single valued map on 
$\Omega_2$, and it can be viewed as a function of $z_1^m z_2^n$.
\[\pi(z_1,z_2):=\] 
\[\left(((z_1^m z_2^n)^{1/(kn+lm)})^l \left(\left(\frac{p_1^k}
{p_2^l}\right)^{1/(kn+lm)}\right)^n,
((z_1^m z_2^n)^{1/(kn+lm)})^k \left(\left(\frac{p_1^k}{p_2^l}\right)^
{1/(kn+lm)}\right)^{-m}\right),\]
where in both of the coordinates the same branch of the root is taken.
\end{Lem}

\begin{proof}
This follows from an easy computation. Since $f_{i,j}$ has a $kn+lm$-symmetry
in the two variables, it is well defined and holomorphic.
\end{proof}

\begin{Lem}{\label{Lem:deelbdd}}
For every $0 \leq i,j \leq kn+lm -1$ there exist functions $f_{i,j}^1$, 
$f_{i,j}^2 \in H^{\infty}(\Omega_2)$ such that \[z_1^i z_2^j f_{i,j}(z_1,z_2)
= f_{i,j}^1(z_1,z_2) \left(\frac{z_1^k}{z_2^l} - \frac{p_1^k}{p_2^l}\right) + 
f_{i,j}^2(z_1,z_2) (z_1^m z_2^n - p_1^m p_2^n).\]
\end{Lem}

\begin{proof}
We start by constructing good holomorphic candidates for $f_{i,j}^1$ and 
$f_{i,j}^2$. Then we show that these functions are indeed bounded.\\
A meromorphic solution of the problem is \[z_1^i z_2^j f_{i,j}(z_1,z_2) = 
z_1^i z_2^j 
\frac{f_{i,j}(z_1,z_2)} {\frac{z_1^k}{z_2^l} - \frac{p_1^k}{p_2^l}} 
\left(\frac{z_1^k}{z_2^l} - \frac{p_1^k}{p_2^l}\right)+ 0 \cdot 
(z_1^m z_2^n - p_1^m p_2^n).\] We search for a function $h$ such that 
\[f_{i,j}^1(z_1,z_2) = z_1^i z_2^j \frac{f_{i,j}(z_1,z_2)}
{\frac{z_1^k}{z_2^l} - \frac{p_1^k}{p_2^l}} + h(z_1,z_2) (z_1^m z_2^n - p_1^m 
p_2^n) \quad \quad \quad \quad \text {(*)}\] and  
\[f_{i,j}^2(z_1,z_2) = - h(z_1,z_2) \left(\frac{z_1^k}{z_2^l} - 
\frac{p_1^k}{p_2^l}\right) \] are holomorphic. Then
\[h(z_1,z_2) = \frac{-f_{i,j}^2(z_1,z_2)}{\frac{z_1^k}{z_2^l} - 
\frac{p_1^k}{p_2^l}},\] \[f_{i,j}^1 (z_1,z_2) = \frac{z_1^i z_2^j f_{i,j}
(z_1,z_2) - f_{i,j}^2(z_1,z_2) (z_1^m z_2^n -p_1^m p_2^n)}
{\frac{z_1^k}{z_2^l} - \frac{p_1^k}{p_2^l}}.\] We want $f_{i,j}^1$ to be 
holomorphic.
Then it is necessary and sufficient that $f_{i,j}^2 (z_1,z_2) = \frac{z_1^i 
z_2^j f_{i,j}(z_1,z_2)}{z_1^m z_2^n - p_1^m p_2^n}$ for points on the zero 
set of 
$\frac{z_1^k}{z_2^l} - \frac{p_1^k}{p_2^l}$. Therefore we define $f_{i,j}^2$ 
as 
\[f_{i,j}^2 (z_1,z_2) := \frac{z_1^i z_2^j f_{i,j} (\pi(z_1,z_2))}{z_1^m z_2^n 
- p_1^m p_2^n},\] and $f_{i,j}^1$ according to (*) as
\[f_{i,j}^1 (z_1,z_2) := \frac{z_1^i z_2^j \left(f_{i,j}(z_1,z_2) - f_{i,j}
(\pi(z_1,z_2))\right)}{\frac{z_1^k}{z_2^l} - \frac{p_1^k}{p_2^l}}.\] 
These are holomorphic functions, and we have that \[z_1^i z_2^j f_{i,j}(z_1,
z_2) = f_{i,j}^1(z_1,z_2) (\frac{z_1^k}{z_2^l} - \frac{p_1^k}{p_2^l}) + 
f_{i,j}^2(z_1,z_2) (z_1^m z_2^n - p_1^m p_2^n).\]
We proceed to show that the functions $f_{i,j}^1$ and $f_{i,j}^2$ are bounded
on $\Omega_2$. 
We start with the function $f_{i,j}^2$.
We define a function $F$, similar to $f_{i,j}^2$, and show that it is bounded
on $\Omega_2$.
\[F(z_1,z_2):= (kn + lm)^2 f_{i,j}^2(z_1,z_2) = \]\[\frac{(\frac{p_1^k z_1^k}
{p_2^l z_2^l})^{\frac{in-jm}{kn+lm}}
\sum_{s,t =1}^{kn + lm} \zeta^{-is -jt} f \left(\zeta^s \left( (z_1^m z_2^n)^l 
(\frac{p_1^k}{p_2^l})^n \right)^{1/(kn+lm)}, \zeta^t \left((z_1^m z_2^n)^k 
(\frac{p_1^k}{p_2^l})^{-m}\right)^{1/(kn+lm)}\right)} {z_1^m z_2^n - p_1^m 
p_2^n}.\] Then $(\frac{p_1^k z_1^k}{p_2^l z_2^l})^{jm - in} F^{kn + lm}$ is 
equal to
\[\left(\frac {\sum_{s,t =1}^{kn + lm} \zeta^{-is -jt} f(\zeta^s \left(
(z_1^m z_2^n)^l (\frac{p_1^k}{p_2^l})^n \right)^{1/(kn+lm)}, \zeta^t 
\left(( z_1^m z_2^n)^k (\frac{p_1^k} {p_2^l})^{-m}\right)^{1/(kn+lm)})}
{z_1^m z_2^n - p_1^m p_2^n}\right)^{kn +lm}.\]
We substitute $x=z_1^n z_2^m$ in the last line, and it becomes
\[\left(\frac {\sum_{s,t =1}^{kn + lm} \zeta^{-is -jt} f(\zeta^s \left(x^l 
(\frac{p_1^k}{p_2^l})^n \right)^{1/(kn+lm)},\zeta^t \left( x^k (\frac{p_1^k}
{p_2^l})^{-m}\right)^{1/(kn+lm)})}{x - p_1^m p_2^n}\right)^{kn +lm}.\]
The numerator is bounded, and we have a removable singularity at $x=p_1^m 
p_2^n$. Hence the 
function is bounded. Since $(\frac{z_1^k}{z_2^l})^{in - jm}$ is bounded, $F$ 
is bounded as well. The same goes for $f_{i,j}^2$.\\
Now we turn our attention to the function $f_{i,j}^1$. Remember that $\omega_2$
was given by $\{(x,y) : y \leq \frac{-m}{n} x + r\}$. The line given by $y = 
\frac{-m}{n} x + r$ corresponds to a curve in $\CC^2$ given by $z_1^m z_2^n = 
nr$. For $|K| \leq nr$, let $\Omega_2^K := \Omega_2 \cap [z_1^m z_2^n =K]$. We 
will estimate $f_{i,j}^1$ on the sets $\Omega_2^K$. Since we have that 
$z_1^i z_2^j f_{i,j}(z_1,z_2)$ is bounded (by construction) and that $z_1^i 
z_2^j f_{i,j} (\pi (z_1,z_2))$ is bounded (as shown while proving that 
$f_{i,j}^2$ is bounded), for every $\mu > 0$ there exists a constant $C$ such 
that 
\[|f_{i,j}^1(z_1,z_2)| \leq C \quad \text { on } \Omega_2 \setminus \{(z_1,z_2)
\in \Omega_2 : \left|\frac{z_1^k}{z_2^l} - \frac{p_1^k}{p_2^l} \right| < 
\mu\}.\] The construction of $\Omega_2$
implies the existence of an $\mu >0$ such that for every $(z_1,z_2) \in 
\Omega_2^K$ (with $|K| \leq nr$), $\Theta \in [0, 2 \pi]$, there is a point
$(s,t) \in \Omega_2^K$ with $\frac{s^k}{t^l} - \frac{p_1^k}{p_2^l} = \mu
e^{i \Theta}$. 
Since $\Omega_2^K$ can locally be seen as an open set in $\CC$ (after the 
appropriate biholomorphic mapping), applying the maximum principle yields that
\[|f_{i,j}^1(z_1,z_2)| \leq \max_{s,t} |f_{i,j}^1(s,t)| \leq C.\]
So $f_{i,j}^1$ is bounded as well.
\end{proof}

\begin{Prop}{\label{Prop:locsol}}
Let $f$ be a bounded holomorphic function on $\Omega_2$ that vanishes at
$(p_1,p_2)$. There exist functions $\tilde{f_1}$, $\tilde{f_2} \in H^{\infty}
(\Omega_2)$ such that \[f(z_1,z_2) = \tilde{f_1} (z_1,z_2) (z_1-p_1) + 
\tilde{f_2} (z_1,z_2) (z_2-p_2).\]
\end{Prop}

\begin{proof}
Since $z_1^m z_2^n - p_1^m p_2^n$ is a polynomial that vanishes at $p$, there
are by lemma \ref{Lem:solpol} polynomials $P_1$ and $P_2$ such that
\[z_1^m z_2^n - p_1^m p_2^n = P_1(z_1,z_2) (z_1-p_1) + P_2(z_1,z_2) (z_2-p_2) 
\quad \quad \forall z_1,z_2 \in \CC^2.\]
Use lemma \ref{Lem:solfrac}
to obtain a similar result for $\frac{z_1^k}{z_2^l} - \frac{p_1^k}{p_2^l}$.
We substitute this and the solutions obtained for $z_1^i z_2^j f_{i,j}(z_1,z_2)
\; (0 \leq i,j \leq kn +lm -1$; note that we may assume that $f_{i,j}(p)=0$,
as remarked after lemma \ref{Lem:symmetry}). This yields that
\[f(z_1,z_2) = \sum_{i,j=0}^{kn +lm -1} z_1^i z_2^jf_{i,j}(z_1,z_2) = \]
\[\sum_{i,j=0}^{kn +lm 
-1} \left( f_{i,j}^1(z_1,z_2) (\frac{z_1^k}{z_2^l} - \frac{p_1^k}{p_2^l}) + 
f_{i,j}^2(z_1,z_2) (z_1^m z_2^n - p_1^m p_2^n) \right) = \]
\[\left(\sum_{i,j=0}^{kn +lm -1} f_{i,j}^1(z_1,z_2)R_1(z_1,z_2) + 
f_{i,j}^2 P_1(z_1,z_2) \right)(z_1-p_1) + \] 
\[\left(\sum_{i,j=0}^{kn +lm -1} f_{i,j}^1(z_1,z_2)R_2(z_1,z_2) + 
f_{i,j}^2 P_2(z_1,z_2) \right) (z_2-p_2) = \]
\[\tilde{f_1} (z_1,z_2) (z_1-p_1) + \tilde{f_2} (z_1,z_2) (z_2-p_2) .\]
The functions $\tilde{f_1}$ and $\tilde{f_2}$ are in $H^{\infty}(\Omega_2)$.
\end{proof}

\section{Main result}

\begin{Thm}{\label{Thm:adomsol}}
Let $\Omega$ be an $A$-domain. Then for every $f \in H^{\infty}(\Omega)$ that 
vanishes at $p=(p_1, p_2) \in \Omega$
there exist functions $f_1$, $f_2 \in H^{\infty}(\Omega)$ such that 
\[f(z_1,z_2) = f_1(z_1,z_2) (z_1-p_1) + f_2(z_1,z_2) (z_2-p_2) \quad \forall z
\in \Omega.\]
Thus one can solve the Gleason problem for $H^{\infty}(\Omega)$.
\end{Thm}

\begin{proof}
Let $\Omega_1$, $\Omega_2$ be as in section \ref{section:pieces}. As noted
there, one can find such $f_1$, $f_2$ if $p \in \Omega_1$. So suppose $p \in
\Omega_2$. We make the local solutions on $\Omega_1$ and $\Omega_2$, using
theorem \ref{Prop:locsol}. The $\overline{\partial}$-problem corresponding to 
the patching of the two local solutions yields a bounded $(0,1)$-form that
has support outside $B(0, \epsilon)$. Theorem \ref{Thm:aholder} yields a
bounded solution to this particular $\overline{\partial}$-problem. Now 
proceed in the standard way (e.g. lemma \ref{Lem:holsupbdd}) to obtain the 
appropriate $f_1$ and $f_2 \in H^{\infty}(\Omega)$.
\end{proof}

\section{The Hartogs triangle and related domains}{\label{S:Hartogsfriends}} 

\noindent
For $k$, $l \in \NN^+$ let $\Omega_{k,l}$ be the domain
defined by \[\Omega_{k,l} := \{(z_1, z_2) \in \CC^2 : 
|z_1|^k < |z_2|^l < 1\}.\] The 
Hartogs triangle is exactly $\Omega_{1,1}$.
The situation becomes slightly more complicated compared to the previous 
sections, since $\Omega_{k,l}$ contains
points of the form $(0,a)$. Thus the functions $f_{i,j}$ as constructed in
lemma \ref{Lem:symmetry} may no longer be holomorphic. We will show that one
can still solve the Gleason problem for $H^{\infty}(\Omega_{k,l})$.\\
If $p_1 \neq 0$, we return to the construction in section 
\ref{section:constrloc}; the domain is now cut off with $|z_2|=1$. We still
project a point of $z_1^m z_2^n = c$ onto the zero set of $\frac{z_1^k}{z_2^l}
- \frac{p_1^k}{p_2^l}$, but now $m=0$, $n=1$, thus $z_1^m z_2^n$ is simply 
$z_2$. Now repeat the proof in section \ref{section:constrloc} to see that
there exist $f_1$, $f_2 \in H^{\infty}(\Omega_{k,l})$
with $f(z) = f_1(z)(z_1 - p_1) + f_2(z)(z_2 - p_2)$.
There are only two things to point out :
\begin{itemize}
\item The functions $f_{i,j}$ may no longer be holomorphic (in their 
definitions we 
divide by $z_1^i$), but $z_1^i z_2^j f_{i,j}$ is still bounded and holomorphic.
\item The expression $\left(\frac{p_1^kz_1^k}{p_2^l z_2^l}\right)^
{\frac{in -jm}{kn + lm}}$ becomes $\left(\frac{p_1^kz_1^k}{p_2^l z_2^l}\right)^
{i/k}$. Thus we only need that $\frac{z_1^k}{z_2^l}$ is bounded, and not that 
$\frac{z_2^l}{z_1^k}$ is bounded. 
\end{itemize}
Now we consider the case that $p_1 = 0$. It is tempting to repeat the previous
argument, but this is impossible. Namely, in the remark after 
\ref{Lem:symmetry}, we assume that $f_{i,j}$ vanishes at $p$. Unfortunately, 
$f_{i,j}$ is not defined at $p$. There is another construction however.

\begin{Lem}
Let $f \in H^{\infty}(\Omega_{k,l})$ such that $f$ 
vanishes at $(0,p_2)$. Let
\[f_1(z_1,z_2) := \frac{z_2^l}{p_2^l} \frac{f(z_1,z_2) - f(0,z_2)}{z_1},\]
\[f_2(z_1,z_2) := \frac{p_2^l - z_2^l}{p_2^l(z_2-p_2)} (f(0,z_2) - f(z_1,z_2)) 
+ \frac{f(0,z_2)}{z_2-p_2}.\]
Then $f_1$, $f_2 \in H^{\infty}(\Omega_{k,l})$ and 
\[f(z_1,z_2) = f_1(z_1,z_2)z_1 + f_2(z_1,z_2) (z_2-p_2) \quad \forall \;
(z_1,z_2) \in \Omega_{k,l}.\]
\end{Lem}

\begin{proof}
We see immediately that $f_1$ and $f_2$ are holomorphic, that $f_2$ is bounded
and that the last equality holds.  We rewrite $f_1$ :
\[f_1 (z_1,z_2) = \frac{z_1^{k - 1}}{p_2^l} \frac{z_2^l}{z_1^{k}}
(f(z_1,z_2) - f(0,z_2)).\]
For every $c \in \CC$ with $|c| \leq 1$
the set $\Omega_{k,l} \cap [z_2=c]$ (a disc) contains a circle with radius 
$|\frac{c}{2}|^{l/k}$. On this circle, we have that 
$|\frac{z_2^l} {z_1^{k}}| = 2^l$. Applying the maximum principle 
yields that
\[|f_1(z_1,c)| \leq  \frac{2^{l+1} ||f||_{\infty}}{p_2^l} \quad \forall \; 
|z_1| \leq |\frac{c}{2}|^{l/k}.\]
It follows that $f_1$ is bounded on $\Omega_{k,l}$.
\end{proof}

\noindent Thus we have the following theorem :

\begin{Thm}{\label{Thm:hartogs}}
For $k$, $l \in \NN^+$ let $\Omega_{k,l}$ be the domain defined by 
\[\Omega_{k,l} := \{(z_1, z_2) \in \CC^2 : |z_1|^k < |z_2|^l < 1\}.\]
One can solve the Gleason problem for $H^{\infty}(\Omega_{k,l})$.
\end{Thm}

\section{If the domain meets one of the coordinate axes}

\noindent In this section, we study domains that are connected both to the 
$A$-domains and the domains $\Omega_{k,l}$. Namely, let $\Omega \subset \CC^2$ 
be a bounded pseudoconvex
Reinhardt domain, such that for $c$ close to $- \infty$, $\partial \omega \cap 
[y<c]$ consists of 2 arcs, one of them being a half line with rational slope. 
We assume that $0 
\notin \Omega$, and that $\Omega$ meets the $z_2$-axis. (Because of symmetry,
everything will hold if $\Omega$ only meets the $z_1$-axis as well.)\\
Let $K_1$, $K_2$ be constants such that $\partial \omega$ is a half line for
$[y < K_1]$, and $\omega$ is rounded off strictly convexily above $[y < K_1]$,
such that $\{(x,y) \in \partial \omega : y > K_1\}$ has $y=K_2$ as horizontal
asymptote. We fix $p = (p_1, p_2) \in \Omega$, and an $f \in H^{\infty}
(\Omega)$ that vanishes at $p$. We will now solve the Gleason problem for $f$
at $p$. There is a strictly convex point $A=(a_1,a_2) \in \partial 
\omega$ with $\log |p_2| < a_2$. This point has a neighborhood in $\partial
\omega$ consisting only of strictly convex points. Take a point $B=(b_1,b_2)$
in this neighborhood with $\log |p_2| < b_2 < a_2$. Then is $\frac{f(z_1,z_2)}
{z_2 - p_2}$ bounded on $\Omega \cap [|z_2| > \exp(b_2)]$, and on this set we 
have that $f(z_1,z_2) = \frac{f(z_1,z_2)}{z_2 - p_2} (z_2-p_2)$. The boundary 
of $\omega \cap [y < a_2]$ is a straight line for $y$ small. Thus one can 
solve 
the Gleason problem for $H^{\infty}(\Omega \cap [|z_2|<\exp(a_2)])$, just as 
section \ref{S:Hartogsfriends}.
One can patch the two local solutions together to a global solution using the
standard techniques, since $\partial \Omega \cap [|z_2| > \exp(b_2)] \cap 
[|z_2|<\exp(a_2)] \subset S(\Omega)$.

\noindent
The case where a part of $\partial \omega$ is described by $[y=c]$ can be dealt
with in a similar way. This yields the following theorem :

\begin{Thm}
Let $\Omega \subset \CC^2$ be a bounded pseudoconvex Reinhardt domain, that
meets exactly one of the axes. Suppose that one part of $\partial \omega$ is 
a half line, and that the other boundary points of $\partial \omega$ are 
strictly convex and $C^5$. Then one can solve the Gleason problem for
$H^{\infty}(\Omega)$.
\end{Thm}

\section{Final remarks}

\noindent The results in this article all rely on theorem \ref{Thm:aholder}. 
As noted before, one can
prove this theorem for Reinhardt domains that for small $z$ look like 
\[\{(z_1,z_2) : a < \left|\frac{z_1^k}{z_2^l}\right| < b\},\]
and are rounded off pseudoconvexily. Thus one can still solve the Gleason
problem if there are ``enough'' strictly pseudoconvex points in the sense of
\cite{lewi}. The condition that the strictly pseudoconvex points
have to be $C^5$ can, as usual, be relaxed to $C^2$, but this would even need
more machinery.

\vskip5mm
\noindent Now let $\Omega$ be a bounded pseudoconvex Reinhardt domain in 
$\CC^2$ that has a strictly pseudoconvex $C^5$ boundary outside a ball around 
the origin. If for $c$ close to $- \infty$, $\partial \omega \cap [y<c]$ 
consists of 2 arcs that have parallel asymptotes  with rational slope, theorem 
\ref{Thm:adomsol} holds for $\Omega$ as well. This is because we are either in
the situation described in the previous remark, or every point in $\omega$ 
lies in a triangle of strictly convex points of $\omega$, and one can apply 
lemma \ref{Lem:driehoek}.

\vskip5mm \noindent
We do not yet know how to solve the Gleason problem for $H^{\infty}(\Omega)$
if $\Omega$ is a Reinhardt domain that for small $z$ looks like 
\[\{(z_1,z_2) : a < \left|\frac{z_1^{\alpha}}{z_2}\right| < b\},\]
where $\alpha \notin \QQ$, or (with $r \not= l$) 
\[\{(z_1,z_2) : a |z_2^l| < |z_1^k| < b |z_2^r|\}.\]
The first problem is hard because $z_1^{\alpha}$ is not a holomorphic function;
the second problem is hard because the function $|z_1^k z_2^{-l}|^{\frac{in 
-jm}{kn +lm}}$ (that appeared in the proof of theorem \ref{Lem:deelbdd}) is
no longer bounded.

\vskip5mm \noindent
Oscar Lemmers

\vskip10mm \noindent
Jan Wiegerinck\\
Department of mathematics\\
University of Amsterdam\\
Plantage Muidergracht 24\\
1018 TV   Amsterdam\\
The Netherlands

\end{document}